\magnification\magstep1
\baselineskip15pt
\def\ref#1{\expandafter\edef\csname#1\endcsname}
%
%
\ref {Introduction}{1}
\ref {Vanishing}{2}
\ref {Fund}{Lemma\penalty 10000\ 2.1}
\ref {IntNsym}{Theorem\penalty 10000\ 2.2}
\ref {IntSym}{Theorem\penalty 10000\ 2.3}
\ref {Constr}{Theorem\penalty 10000\ 2.4}
\ref {lam*}{Corollary\penalty 10000\ 2.5}
\ref {SymNonsym}{Corollary\penalty 10000\ 2.6}
\ref {Hecke}{3}
\ref {Van}{Lemma\penalty 10000\ 3.1}
\ref {Ei}{Theorem\penalty 10000\ 3.6}
\ref {Extra vanishing}{4}
\ref {L1}{Lemma\penalty 10000\ 4.1}
\ref {Stab}{Theorem\penalty 10000\ 4.4}
\ref {T1}{Theorem\penalty 10000\ 4.5}
\ref {T2}{Theorem\penalty 10000\ 4.6}
\ref {Inversion}{5}
\ref {Diag}{Theorem\penalty 10000\ 5.1}
\ref {Classical}{6}
\ref {References}{7}
\ref {Ch}{[Ch]}
\ref {HU}{[HU]}
\ref {Kn1}{[Kn]}
\ref {KS}{[KS]}
\ref {M1}{[M1]}
\ref {M2}{[M2]}
\ref {Ok}{[Ok]}
\ref {Ol}{[Ol]}
\ref {Op}{[Op]}
\ref {Sa}{[Sa]}
\ref {Wa}{[WUN]}
\def\today{\number\day.~\ifcase\month\or
  Januar\or Februar\or M\"arz\or April\or Mai\or Juni\or
  Juli\or August\or September\or Oktober\or November\or Dezember\fi
  \space\number\year}
\font\sevenex=cmex7
\scriptfont3=\sevenex
\font\fiveex=cmex10 scaled 500
\scriptscriptfont3=\fiveex

\def\ZS{{\widetilde Z}}
\def\phi{\varphi}
\def\epsilon{\varepsilon}
\def\theta{\vartheta}

\def\uauf{\lower1.7pt\hbox to 3pt{%
\vbox{\offinterlineskip
\hbox{\vbox to 8.5pt{\leaders\vrule width0.2pt\vfill}%
\kern-.3pt\hbox{\lams\char"76}\kern-0.3pt%
$\raise1pt\hbox{\lams\char"76}$}}\hfil}}
\def\cite#1{\expandafter\ifx\csname#1\endcsname\relax
{\bf?}\immediate\write16{#1 ist nicht definiert!}\else\csname#1\endcsname\fi}
\def\expandwrite#1#2{\edef\next{\write#1{#2}}\next}
\def\neverexpand{\noexpand\noexpand\noexpand}
\def\strip#1\ {}
\def\ncite#1{\expandafter\ifx\csname#1\endcsname\relax
{\bf?}\immediate\write16{#1 ist nicht definiert!}\else
\expandafter\expandafter\expandafter\strip\csname#1\endcsname\fi}
\newwrite\AUX
\immediate\openout\AUX=\jobname.aux
\newcount\Abschnitt\Abschnitt0
\def\beginsection#1. #2 \par{\advance\Abschnitt1%
\vskip0pt plus.10\vsize\penalty-250
\vskip0pt plus-.10\vsize\bigskip\vskip\parskip
\edef\TEST{\number\Abschnitt}
\expandafter\ifx\csname#1\endcsname\TEST\relax\else
\immediate\write16{#1 hat sich geaendert!}\fi
\expandwrite\AUX{\neverexpand\ref{#1}{\TEST}}
\leftline{\bf\number\Abschnitt. \ignorespaces#2}%
\nobreak\smallskip\noindent\SATZ1}
\def\Proof:{\par\noindent{\it Proof:}}
\def\Remark:{\ifdim\lastskip<\medskipamount\removelastskip\medskip\fi
\noindent{\bf Remark:}}
\def\Remarks:{\ifdim\lastskip<\medskipamount\removelastskip\medskip\fi
\noindent{\bf Remarks:}}
\def\Definition:{\ifdim\lastskip<\medskipamount\removelastskip\medskip\fi
\noindent{\bf Definition:}}
\def\Example:{\ifdim\lastskip<\medskipamount\removelastskip\medskip\fi
\noindent{\bf Example:}}
\newcount\SATZ\SATZ1
\def\proclaim #1. #2\par{\ifdim\lastskip<\medskipamount\removelastskip
\medskip\fi
\noindent{\bf#1.\ }{\it#2}\Par
\ifdim\lastskip<\medskipamount\removelastskip\goodbreak\medskip\fi}
\def\Aussage#1{%
\expandafter\def\csname#1\endcsname##1.{\ifx?##1?\relax\else
\edef\TEST{#1\penalty10000\ \number\Abschnitt.\number\SATZ}
\expandafter\ifx\csname##1\endcsname\TEST\relax\else
\immediate\write16{##1 hat sich geaendert!}\fi
\expandwrite\AUX{\neverexpand\ref{##1}{\TEST}}\fi
\proclaim {\number\Abschnitt.\number\SATZ. #1\global\advance\SATZ1}.}}
\Aussage{Theorem}
\Aussage{Proposition}
\Aussage{Corollary}
\Aussage{Lemma}
\font\la=lasy10
\def\strich{\hbox{$\vcenter{\hbox
to 1pt{\leaders\hrule height -0,2pt depth 0,6pt\hfil}}$}}
\def\dashedrightarrow{\hbox{%
\hbox to 0,5cm{\leaders\hbox to 2pt{\hfil\strich\hfil}\hfil}%
\kern-2pt\hbox{\la\char\string"29}}}

\def\Bindestrich{\penalty10000-\hskip0pt}
\let\_=\Bindestrich
\def\.{{\sfcode`.=1000.}}

\def\Rechts#1{\rlap{$\scriptstyle#1$}}
\def\Par{\par}
\def\:={\mathrel{\raise0,9pt\hbox{.}\kern-2,77779pt
\raise3pt\hbox{.}\kern-2,5pt=}}
\def\=:{\mathrel{=\kern-2,5pt\raise0,9pt\hbox{.}\kern-2,77779pt
\raise3pt\hbox{.}}} 

\def\pfeil{\rightarrow}

\def\pf#1{\buildrel#1\over\rightarrow}
\def\Pf#1{\buildrel#1\over\longrightarrow}

\def\Ugleich{\hbox{$\cup$\kern.5pt\vrule depth -0.5pt}}
\def\|#1|{\mathop{\rm#1}\nolimits}
\def\<{\langle}
\def\>{\rangle}
\let\Times=\times
\def\times{\mathop{\Times}}
\let\Otimes=\otimes
\def\otimes{\mathop{\Otimes}}
\catcode`\@=11
\def\hex#1{\ifcase#1 0\or1\or2\or3\or4\or5\or6\or7\or8\or9\or A\or B\or
C\or D\or E\or F\else\message{Warnung: Setze hex#1=0}0\fi}
\def\fontdef#1:#2,#3,#4.{%
\alloc@8\fam\chardef\sixt@@n\FAM
\ifx!#2!\else\expandafter\font\csname text#1\endcsname=#2
\textfont\the\FAM=\csname text#1\endcsname\fi
\ifx!#3!\else\expandafter\font\csname script#1\endcsname=#3
\scriptfont\the\FAM=\csname script#1\endcsname\fi
\ifx!#4!\else\expandafter\font\csname scriptscript#1\endcsname=#4
\scriptscriptfont\the\FAM=\csname scriptscript#1\endcsname\fi
\expandafter\edef\csname #1\endcsname{\fam\the\FAM\csname text#1\endcsname}
\expandafter\edef\csname hex#1fam\endcsname{\hex\FAM}}
\catcode`\@=12 

\fontdef Ss:cmss10,,.
\fontdef Fr:eufm10,eufm7,eufm5.


\def\fQ{{\Fr Q}}

\def\fs{{\Fr s}}

\newread\AUXX
\immediate\openin\AUXX=msxym.tex
\ifeof\AUXX
\fontdef bbb:msbm10,msbm7,msbm5.
\fontdef mbf:cmmib10,cmmib7,.
\else
\fontdef bbb:msym10,msym7,msym5.
\fontdef mbf:cmmib10,cmmib10 scaled 700,.
\fi
\immediate\closein\AUXX

\def\NN{{\bbb N}}
\def\QQ{{\bbb Q}}

\def\ZZ{{\bbb Z}}
\def\cA{{\cal A}}
\def\cE{{\cal E}}\def\cH{{\cal H}}

\def\cP{{\cal P}}

\mathchardef\leer=\string"0\hexbbbfam3F
\mathchardef\subsetneq=\string"3\hexbbbfam24
\mathchardef\semidir=\string"2\hexbbbfam6E
\mathchardef\dirsemi=\string"2\hexbbbfam6F
\let\OL=\overline
\def\overline#1{{\hskip1pt\OL{\hskip-1pt#1\hskip-1pt}\hskip1pt}}
\def\Aq{{\overline{A}}}

\def\Eq{{\overline{E}}}
\def\fQ{{\overline{f}}}

\def\Hq{{\overline{H}}}

\def\Pq{{\overline{P}}}

%
\abovedisplayskip 9.0pt plus 3.0pt minus 3.0pt
\belowdisplayskip 9.0pt plus 3.0pt minus 3.0pt
\newdimen\Grenze\Grenze2\parindent\advance\Grenze1em
\newdimen\Breite
\newbox\DpBox
\def\NewDisplay#1$${\Breite\hsize\advance\Breite-\hangindent
\setbox\DpBox=\hbox{\hskip2\parindent$\displaystyle{#1}$}%
\ifnum\predisplaysize<\Grenze\abovedisplayskip\abovedisplayshortskip
\belowdisplayskip\belowdisplayshortskip\fi
\global\futurelet\nexttok\WEITER}
\def\WEITER{\ifx\nexttok\qed\expandafter\leftQEDdisplay
\else\leftdisplay\fi}
\def\leftdisplay{\hskip-\hangindent\leftline{\box\DpBox}$$}
\def\leftQEDdisplay{\hskip-\hangindent
\line{\copy\DpBox\hfill\lower\dp\DpBox\copy\QEDbox}%
\belowdisplayskip0pt$$\bigskip\let\nexttok=}
\everydisplay{\NewDisplay}
\newbox\QEDbox
\newbox\nichts\setbox\nichts=\vbox{}\wd\nichts=2mm\ht\nichts=2mm
\setbox\QEDbox=\hbox{\vrule\vbox{\hrule\copy\nichts\hrule}\vrule}
\def\qed{\leavevmode\unskip\hfil\null\nobreak\hfill\copy\QEDbox\medbreak}
\newdimen\HIindent
\newbox\HIbox
\def\setHI#1{\setbox\HIbox=\hbox{#1}\HIindent=\wd\HIbox}
\def\HI#1{\par\hangindent\HIindent\hangafter=0\noindent\leavevmode
\llap{\hbox to\HIindent{#1\hfil}}\ignorespaces}
\def\rho{\varrho}

\def\lamq{{\overline\lambda}}
\def\muq{{\overline\mu}}
\fontdef ams: msam10,,.
\mathchardef\cle"1\hexamsfam34
\fontdef Ss:cmss10,,.
\font\BF=cmbx10 scaled \magstep2
\font\CSC=cmcsc10 
\baselineskip15pt
{\baselineskip1.5\baselineskip\rightskip0pt plus 5truecm
\leavevmode\vskip0truecm\noindent
\BF Symmetric and Non\_Symmetric Quantum Capelli Polynomials

}
\vskip1truecm
\leftline{{\CSC Friedrich Knop}%
\footnote*{\rm Partially supported by a grant of the NSF}}
\leftline{Department of Mathematics, Rutgers University, New Brunswick NJ
08903, USA}
\leftline{knop@math.rutgers.edu}
\vskip1truecm
\beginsection Introduction. Introduction

Generalizing the classical Capelli identity has recently attracted a lot
of interest (\cite{HU}, \cite{Ok}, \cite{Ol}, \cite{Sa}, \cite{Wa}). In
several of these papers it was realized, in various degrees of
generality, that Capelli identities are connected with certain symmetric
polynomials which are characterized by their vanishing at certain points.
From this point of view, these polynomials have been constructed by Sahi
\cite{Sa} and were studied in \cite{KS}.

The purpose of this paper is twofold: we quantize the vanishing
condition in a rather straightforward manner and obtain a family
of symmetric polynomials which is indexed by partitions and which
depends on two parameters $q$, $t$. As in \cite{KS}, their main feature
is that they are non\_homogeneous and one of our principal results
states that the top degree terms are the Macdonald polynomials. It is
an interesting problem whether these quantized Capelli polynomials are
indeed connected with quantized Capelli identities (see \cite{Wa}) as
it is in the classical case.

But the main progress over \cite{KS} is the introduction of a family of
{non\_symmetric\/} polynomials which are also defined by vanishing
conditions. They are non\_homogeneous and their top degree terms turns
out to be the non\_symmetric Macdonald polynomials. To prove this, we
introduce certain difference operators of Cherednik type of which our
polynomials are a simultaneous eigenbasis. Because of these operators,
the non\_symmetric functions are much easier to handle than the
symmetric ones. Moreover, the latter can be obtained by a simple
symmetrization process.

More specifically, the non\_symmetric vanishing conditions are as
follows: For $\lambda\in\Lambda:=\NN^n$ let
$|\lambda|:=\sum\lambda_i$ and let $w_\lambda$ be the shortest
permutation such that $w_\lambda^{-1}(\lambda)$ is a partition (i.e., a
non\_increasing sequence). Let $q$ and $t$ be two formal parameters and
consider the vector
$\rho:=(1,t^{-1},t^{-2},\ldots,t^{{-}n{+}1})$. Then we prove that for
every $\lambda\in\Lambda$ there is a polynomial
$E_\lambda(z_1,\ldots,z_n)$, unique up to a scalar factor, which is of
degree $|\lambda|$ and which satisfies the following condition:
$$
E_\lambda(q^\mu w_\mu(\rho))=0\quad\hbox{ for all
$\mu\in\Lambda$ with $|\mu|\le|\lambda|$ and $\mu\ne\lambda$.}
$$
We show that the affine Hecke algebra acts on these polynomials in a
natural way and that there are Cherednik\_type operators of which they
are simultaneous eigenfunctions. This gives the link to the
theory of homogeneous (symmetric or not) Macdonald polynomials. 

Furthermore, in this paper we show two results which are in a way dual to
each other: the polynomial $E_\lambda$ contains, in general, much fewer
monomials than there are of degree $|\lambda|$ (triangularity) and it
vanishes at many more points than required by definition (extra
vanishing). The extra vanishing is expressed in terms of an order
relation on $\Lambda$ which generalizes the order of partitions by
inclusion of diagrams.

Later we prove that the quantized Capelli polynomials can be expressed in
terms of their highest degree component (a Macdonald polynomial) and
certain difference operators (inversion formula). This is used to transfer
previous integrality results of mine, \cite{Kn1}, to the case of Capelli
polynomials.

In the final section, we discuss the transition from the quantum to the
classical case. For this we put $t=q^r$ and let $q$ tend to one. 
\medskip\noindent
{\bf Acknowlegment:} I would like to thank G.~Heckman whose questions
initiated the research to this paper and S.~Sahi for many discussions on
the case of the classical limit.

\beginsection Vanishing. The vanishing condition

Let $\Lambda:=\NN^n$ and $\Lambda^+\subseteq\Lambda$ the subset of
partitions. For $\lambda=(\lambda_i)\in\Lambda$ we write
$|\lambda|:=\sum_i\lambda_i$ and
$l(\lambda):=\|max|\{i\mid\lambda_i\ne0\}$ (with $l(0)=0$).

Let $k$ be a field of characteristic zero and $\cP:=k[z_1,\ldots,z_n]$ the
polynomial ring and $\cP':=k[z_1,z_1^{-1},\ldots,z_n,z_n^{-1}]$ the
ring of Laurent polynomials. To each $\lambda\in\Lambda$ corresponds a
monomial $z^\lambda=\prod_iz_i^{\lambda_i}$. Fix two non\_zero
elements $q$ and $t$ of $k$. Throughout the paper we assume that
$q^at^b\ne1$ for all integers $a,b\ge1$.

The symmetric group $W:=S_n$ acts in the obvious way on $\Lambda$ and
$\cP$. Every $\lambda\in\Lambda$ contains a unique partition
$\lambda^+$ in its $W$\_orbit. Let $w_\lambda\in W$ be the shortest
permutation such that $w_\lambda(\lambda^+)=\lambda$. For all
$\lambda\in\Lambda$ and $x=(x_i), y=(y_i)\in k^n$ let
$q^\lambda:=(q^{\lambda_i})$ and $xy:=(x_iy_i)$. Consider the element
$\rho:=(1,t^{-1},t^{-2},\ldots,t^{-n+1})\in k^n$.

For every $\lambda\in\Lambda$ we define
$\lamq:=w_\lambda(q^{\lambda^+}\rho)$. More concretely,
$\lamq_i=q^{\lambda_i}t^{-k_i}$ where
$$
k_i=k_i(\lambda)=\#\{j=1,\ldots,i-1\mid\lambda_j\ge\lambda_i\}+
\#\{j=i+1,\ldots,n\mid\lambda_j>\lambda_i\}.
$$
The following simple lemma is fundamental:

\Lemma Fund. For $\lambda\in\Lambda$ with $\lambda_n\ne0$ let
$\lambda^*\:=(\lambda_n-1,\lambda_1,\ldots,\lambda_{n-1})$. Then
$\lamq^*=(\lamq_n/q,\lamq_1,\ldots,\lamq_{n-1})$.

\Proof: Follows easily from the definition.\qed

\Theorem IntNsym. For $d\in\NN$ let $S(n,d)$ be the set of all $\lamq$
where $\lambda\in\Lambda$ and $|\lambda|\le d$. Let $\fQ:S(n,d)\pfeil k$
be a mapping. Then there exists a unique polynomial $f\in\cP$ of degree
at most $d$ such that $f(z)=\fQ(z)$ for all $z\in S(n,d)$.

\Proof: The cardinality of $S(n,d)$ equals the dimension of the space of
polynomials of degree at most $d$. Hence existence of $f$ will imply its
uniqueness.

To show existence, we proceed by induction on $n+d$. Every polynomial
can be uniquely written as
$$
f(z_1,\ldots,z_n)=g(z_1,\ldots,z_{n-1})+
(z_n-t^{-n+1})h(z_n/q,z_1,z_2,\ldots,z_{n-1}).
$$
Consider first the set $S_0$ of $\lamq\in S(n,d)$ with $\lambda_n=0$,
i.e., $\lamq_n=t^{-n+1}$. Then, as $z$ runs through $S_0$,
$z':=(z_1,\ldots,z_{n-1})$ will run through $S(n-1,d)$. By induction,
one can choose $g$ such that $f$ takes the required values at $S_0$.

Consider now the set $S_1$ of remaining points $\lamq$ with
$\lambda_n\ne0$. By the lemma, as $z$ runs through $S_1$,
$(z_n/q,z_1,z_2,\ldots,z_{n-1})$ will run through $S(n,d-1)$. The factor
$z_n-t^{-n+1}=\lamq_n-t^{-n+1}$ is not zero by the choice of $q$ and
$t$. By induction, we can find $h$ of degree at most $d-1$ with arbitrary
values at $S(n,d-1)$. So $f$ exists.\qed

There is also a statement for symmetric polynomials which is the
quantized version of a theorem of Sahi \cite{Sa}. For
$\lambda\in\Lambda^+$ let $m_\lambda$ be the corresponding monomial
symmetric polynomial.

\Theorem IntSym. For $d\in\NN$ let $S^+(n,d)$ be the set of all
$\lamq=q^\lambda\rho$ where $\lambda\in\Lambda^+$ and $|\lambda|\le
d$. Let $\fQ:S^+(n,d)\pfeil k$ be a mapping. Then there exists a unique
polynomial $f\in\cP$ of degree at most $d$ such that $f(z)=\fQ(z)$ for
all $z\in S^+(n,d)$.

\Proof: The proof is completely analogous to that in \cite{Sa}. Again,
only existence has to be proved. Let $g\mapsto g^+$ be the linear map
from symmetric polynomials in $n-1$ variables to those with $n$ variables
which sends $m_\lambda(z_1-t^{-n+1},\ldots,z_{n-1}-t^{-n+1})$ to 
$m_{\lambda,0}(z_1-t^{-n+1},\ldots,z_n-t^{-n+1})$. This map preserves
degrees and satisfies
$g^+(z_1,\ldots,z_{n-1},t^{-n+1})=g(z_1,\ldots,z_{n-1})$. We construct
$f$ by setting
$$
f(z_1,\ldots,z_n)=g^+(z_1,\ldots,z_n)+
\prod_{i=1}^n(z_i-t^{-n+1})h(z_1/q,\ldots,z_n/q),
$$
where $g^+$ and $h$ are symmetric. Then, as above, evaluation at
$\lamq\in S^+(n,d)$ with $\lambda_n=0$ gives $g$. Evaluation at the
other points gives $h$.\qed

We obtain the following Theorem/Definition:

\Theorem Constr. a) For every $\lambda\in\Lambda$ there is a unique
polynomial $E_\lambda$ with $E_\lambda(\muq)=0$ for all $\mu\in\Lambda$
with $|\mu|\le|\lambda|$, $\mu\ne\lambda$ and which has an expansion
$E_\lambda=\sum_\mu e_{\lambda\mu}z^\mu$ with $e_{\lambda\lambda}=1$.
\Par\noindent
b) For every $\lambda\in\Lambda^+$ there is a unique
symmetric polynomial $P_\lambda$ with $P_\lambda(\muq)=0$ for all
$\mu\in\Lambda^+$ with $|\mu|\le|\lambda|$, $\mu\ne\lambda$ and which
has an expansion $P_\lambda=\sum_\mu p_{\lambda\mu}m_\mu$ with
$p_{\lambda\lambda}=1$.

\Proof: By \cite{IntNsym}, there is a polynomial $E_\lambda$
satisfying the vanishing condition with $E_\lambda(\lamq)\ne0$. We
have to show that it contains $z^\lambda$ with a non\_zero coefficient.
Let
$$
E_\lambda(z_1,\ldots,z_n)=g(z_1,\ldots,z_{n-1})+
(z_n-t^{-n+1})h(z_n/q,z_1,z_2,\ldots,z_{n-1}).
$$
If $\lambda_n=0$, then $g=E_{\lambda'}$ with
$\lambda'=(\lambda_1,\ldots,\lambda_{n-1})$. By induction, $g$, and
therefore $f$, contain $z^\lambda$. If $\lambda_n\ne0$ then $g=0$ and
$h=E_{\lambda^*}$. We conclude again by induction.

The proof in the symmetric case is analogous.\qed

\noindent Our proof actually gives a bit more. Consider the operators
$$
\Delta f(z_1,\ldots,z_n):=f(z_n/q,z_1,\ldots,z_{n-1}),
$$
and $\Phi:=(z_n-t^{-n+1})\Delta$. Then we have

\Corollary lam*. For $\lambda\in\Lambda$ with $\lambda_n\ne0$ let
$\lambda^*=(\lambda_n-1,\lambda_1,\ldots,\lambda_{n-1})$. Then
$E_\lambda=\Phi(E_{\lambda^*})$

For $\lambda\in\Lambda$ let $\cP_\lambda\subseteq\cP$ be the set of
polynomials of degree at most $|\lambda|$ which vanish at all $\muq$
with $|\mu|\le|\lambda|$ and $\mu\not\in W\lambda$. Of course,
$\cP_\lambda$ depends only on $\Lambda^+$.

\Corollary SymNonsym. We have
$\cP=\oplus_{\lambda\in\Lambda^+}\cP_\lambda$. The set
$\{E_{w\lambda}\mid w\in W\}$ forms a basis of $\cP_\lambda$. Moreover,
$\cP_\lambda\cap\cP^W=kP_\lambda$.

\Proof: Follows directly from the definition.\qed

\noindent We conclude this section by giving two examples. For $k\in\NN$
we define the $q$\_factorial polynomial as
$[z;k]_q:=(z-1)(z-q)\ldots(z-q^{k-1})$. Then the following is obvious.

\Proposition. Let $t=1$. Then $E_\lambda(z;q,1)=[z_1;\lambda_1]_q\ldots
[z_n;\lambda_n]_q$ and $P_\lambda(z;q,1)$ is the symmetrization of it.

\noindent Now we consider the case $t=q$ in the symmetric case. For
$\lambda\in\Lambda^+$ we define the {\it $q$\_factorial Schur
function\/} as $\fs_\lambda(z;q):=a^{-1}\|det|[z_i;\lambda_j+n-j]_q$
where $a=\prod_{i<j}(z_i-z_j)$ is the Vandermonde determinant.

\Proposition. Let $t=q$. Then
$P_\lambda(z;q,q)=q^{-(n-1)|\lambda|}\fs_\lambda(q^{n-1}z;q)$.

\Proof: The proof is the same as in the classical case
\cite{KS}~Prop.3.3.\qed

\beginsection Hecke. Hecke operators

In this section, we are constructing operators which are adapted to the
decomposition $\cP=\oplus_{\lambda\in\Lambda^+}\cP_\lambda$. Let
$s_i\in W$ be the $i$-th simple reflection. Then
$$
N_i:={1-s_i\over z_i-z_{i{+}1}}
$$
is a well defined operator on $\cP$. We define the Hecke operators
$$
\eqalign{
H_i&:=s_i-(1-t)N_iz_i=ts_i-(1-t)z_iN_i,\cr
\Hq_i&:=s_i-(1-t)z_{i+1}N_i=ts_i-(1-t)N_iz_{i+1}.\cr}
$$
They satisfy the relations $H_i-\Hq_i=t-1$ and $H_i\Hq_i=t$. In
particular, both $H_i$ and $-\Hq_i$ satisfy the equation
$(x+1)(x-t)=0$. In addition the braid relations hold
$$
\eqalign{
H_iH_{i{+}1}H_i&=H_{i{+}1}H_iH_{i{+}1}\qquad i=1,\ldots, n-2\cr
H_iH_j&=H_jH_i\qquad\qquad|i-j|>1\cr}
$$
This means that the algebra $\cH$ generated by the $H_i$ is a Hecke
algebra of type $A_{n-1}$. For details see \cite{Ch}, \cite{M2}, or
\cite{Kn1}.

\Lemma Van. Let $\mu\in\Lambda$ and $f\in\cP'$.
Then $H_if(\muq)$ and $\Hq_if(\muq)$ are linear combinations of
$f(\muq)$ and $f(\overline{s_i\mu})$ where the coefficients are
independent of $f$.

\Proof: We have
$$
\Hq_if(\muq)=[s_i-(1-t)z_{i+1}N_i]f(\muq)=
{(t-1)\muq_{i+1}\over \muq_i-\muq_{i+1}}f(\muq)+
{\muq_i-t\muq_{i+1}\over \muq_i-\muq_{i+1}}f(s_i\muq).
$$
If $\mu_i\ne\mu_{i+1}$ then $s_i\muq=\overline{s_i\mu}$. Otherwise,
$\muq_i-t\muq_{i+1}=0$ by definition of $\muq$.\qed

\Corollary. For every $\lambda\in\Lambda^+$ holds
$\cH\cP_\lambda\subseteq\cP_\lambda$.

\Corollary. Let $X$ be an operator in the algebra $\cA$ generated by the
$H_i$, $z_i$, and $\Phi$ (respectively in the both sided ideal
$\cA\Phi\cA$). Let $\lambda\in\Lambda$ and $f\in\cP'$. Then $Xf(\lamq)$
is a linear combination of $f(\muq)$ where $\mu\in\Lambda$ and
$|\mu|\le|\lambda|$ (respectively $|\mu|<|\lambda|$).

\Corollary. Let $\lambda\in\Lambda$ with $\lambda_i=\lambda_{i+1}$.
Then $\Hq_i(E_\lambda)=E_\lambda$ and $H_i(E_\lambda)=tE_\lambda$.

\Proof: Let us consider $E:=H_i(E_\lambda)$ and $\mu\in\Lambda$ with
$|\mu|\le|\lambda|$ and $\mu\ne\lambda$. Then
$E_\lambda(\muq)=E_\lambda(\overline{s_i\mu})=0$, hence $E(\muq)=0$.
This means that $E$ is a multiple of $E_\lambda$. Evaluation at
$z=\lamq$ implies that the factor is $t$.\qed

\noindent
For $i=1,\ldots,n$ we define the Cherednik operators
$$
\xi_i^{-1}=\Hq_i\ldots\Hq_{n-1}\Delta H_1\ldots H_{i-1}.
$$
Their relevance will become clear later. Furthermore, we define the
operators
$$
\Xi_i:=z_i^{-1}+z_i^{-1}H_i\ldots H_{n-1}\Phi H_1\ldots H_{i-1}
$$
which are a priori only well\_defined on Laurent polynomials.
The following relations are easily established
$$
H_i\Xi_i=\Xi_{i+1}\Hq_i,\quad i=1,\ldots,n-1;\qquad
H_i\Xi_j=\Xi_jH_i,\quad j\ne i,i+1.
$$

\Lemma. The operators $\Xi_i$ act on $\cP$.

\Proof: Since $\Xi_i=t^{-1}\Hq_i\Xi_{i+1}\Hq_i$ it suffices to
consider the case $i=n$. Because
$$
\Xi_n:=(z_n-t^{-n+1})z_n^{-1}(\xi_n^{-1}-t^{n-1})+t^{n-1},
$$
the assertion follows from the following claim
$$
\xi_n^{-1}f=t^{n-1}f\ \|mod|z_n\quad\hbox{for all $f\in\cP$}
$$
We write $f\equiv g$ for $f=g\ \|mod|z_n$. Then we show first by
induction
$$
\llap{\hbox to 2\parindent{$(*)$\hfill}}
\Delta H_1\ldots H_if\equiv t^{i-1}\Delta s_1\ldots s_i.
$$
We have $\Delta H_1\ldots H_if\equiv t^{i-2}\Delta s_1\ldots
s_{i-1}(ts_i-(1-t)z_iN_i)f$. The lemma follows from $\Delta s_1\ldots
s_{i-1}z_i=q^{-1}z_n\Delta s_1\ldots s_{i-1}$.

For $i=n-1$, $(*)$ reads $\xi_n^{-1}f\equiv t^{n-1}\Delta s_1\ldots
s_{n-1}f$.  The claim follows from
$$
\Delta s_1\ldots s_{n-1}f(z)=f(z_1,\ldots,z_{n-1},z_n/q)\equiv f(z).
$$\qed

\noindent
The $\Xi_i$ are inhomogeneous versions of the Cherednik operators.
Observe, that they don't increase the degree. The main result of this
paper is

\Theorem Ei. For all $\lambda\in\Lambda$ and $i=1,\ldots,n$ holds
$\Xi_i(E_\lambda)=\lamq_i^{-1}E_\lambda$.

\Proof: Let $d:=|\lambda|$ and write $\Xi_i=z_i^{-1}+X_i$.
Then $\Xi_i(E_\lambda)=z_i^{-1}E_\lambda+X_i(E_\lambda)$. Since
$E_\lambda$ vanishes in $S(n,d-1)$, \cite{Van}
implies that $X_i(E_\lambda)$ vanishes in
$S(n,d)$. Hence, $\Xi_i(E_\lambda)$ vanishes for all $\mu\in
S(n,d)$ with $\mu\ne\lambda$. This implies that
$\Xi_i(E_\lambda)=cE_\lambda$ for some $c\in k$. Evaluation at
$z=\lamq$ implies $c=\lamq_i^{-1}$.\qed

\Corollary. The operators $\Xi_1,\ldots,\Xi_n$ commute pairwise.

\Corollary. Let $p\in\cP^W$. Then $\Xi_p:=p(\Xi_1,\ldots,\Xi_n)$
commutes with all $H_i$. Moreover,
$\Xi_p(P_\lambda)=p(\lamq^{-1})P_\lambda$.

\Proof: $\Xi_p$ acts on $\cP_\lambda$ as scalar multiplication by
$p(\lamq^{-1})$.\qed

\Theorem. Let $\Eq_\lambda$ be the top homogeneous part of
$E_\lambda$. Then $\xi_i^{-1}(\Eq_\lambda)=\lamq_i^{-1}\Eq_\lambda$.
This means, that $\Eq_\lambda$ is a non\_symmetric Macdonald
polynomial. Similarly, the top homogeneous part $\Pq_\lambda$ of
$P_\lambda$ is a symmetric Macdonald polynomial.

\Proof: We have $\Xi_n=\xi_n^{-1}$ plus terms which decrease the
degree. Therefore, also $\Xi_i=\xi_i^{-1}$ plus degree decreasing
operators. The theorem follows since the Macdonald polynomials are
characterized as eigenfunctions of the $\xi_i$.\qed

There is a (partial) order relation on
$\Lambda$. First, recall the usual order on the set
$\Lambda^+$: we say
$\lambda\ge\mu$ if $|\lambda|=|\mu|$ and
$$
\lambda_1+\lambda_2+\ldots+\lambda_i\ge\mu_1+\mu_2+\ldots+\mu_i
\quad\hbox{for all }i=1,\ldots,n.
$$
This order relation is extended to all of $\Lambda$ as follows.
For every $\lambda\in\Lambda$ there is a unique partition
$\lambda^+$ in the orbit $W\lambda$. For all permutations $w\in W$ with
$\lambda=w\lambda^+$ there is a unique one, denoted by $w_\lambda$, of
minimal length. We define $\lambda\ge\mu$ if either $\lambda^+>\mu^+$
or $\lambda^+=\mu^+$ and $w_\lambda\le w_\mu$ in the Bruhat order of
$W$. In particular, $\lambda^+$ is the unique {\it maximum\/} of
$W\lambda$.

\Lemma. The operators $\Xi_i$ are triangular. More precisely,
$\Xi_i(z^\lambda)=
\lamq_i^{-1}z^\lambda+\sum_{\mu<\lambda}c_{\lambda\mu}z^\mu$.

\Proof: For this we write $\Xi_i=z_i^{-1}+\Hq_i\ldots\Hq_{n-1}
(1-t^{-n+1}z_n^{-1})\Delta H_1\ldots H_{i-1}=\xi_i^{-1}+Y_i$ where
$Y_i=z_i^{-1}-t^{-n+1}\Hq_i\ldots\Hq_{n-1}z_n^{-1}\Delta H_1\ldots
H_{i-1}$. It is well known that $\xi_i^{-1}$ is triangular (see
\cite{M2}) with the given coefficient of $z^\lambda$. Since $Y_i$
decreases the degree it suffices to show that $\mu^+<\lambda^+$ for
each monomial $z^\mu$ occuring in $Y_i(z^\lambda)$. But that is
well known to be true.\qed

\Theorem. For every $\lambda\in\Lambda$ there is the expansion
$E_\lambda=z^\lambda+\sum_{\mu\in\Lambda\atop\mu<\lambda}
e_{\lambda\mu}z^\mu$.
Similarly, if $\lambda\in\Lambda^+$ then $P_\lambda=
m_\lambda+\sum_{\mu\in\Lambda^+\atop\mu<\lambda\hfill}
p_{\lambda\mu}m_\mu$.

\Proof: By the triangularity and diagonalizability of $\Xi_i$ there must
be an eigenfunction of the stated form with eigenvalue $\lamq_i^{-1}$.
Thus, it equals $E_\lambda$.\qed

\beginsection Extra vanishing. The extra vanishing theorem

We are going to introduce another (partial) order relation on $\Lambda$.
Let $\lambda,\mu\in\Lambda$. Then we say $\lambda\cle\mu$ if there is
a permutation $\pi\in W$ such that $\lambda_i<\mu_{\pi(i)}$ if
$i<\pi(i)$ and $\lambda_i\le\mu_{\pi(i)}$ if $i\ge\pi(i)$. In this case
we call $\pi$ a {\it defining permutation\/} for $\lambda\cle\mu$.

\Lemma L1. If $\lambda\cle\mu$ and $|\lambda|\ge|\mu|$ then
$\lambda=\mu$. 

\Proof: All inequalities $\lambda_i\le\mu_{\pi(i)}$ must the equalities.
This can only happen if $i\ge\pi(i)$ for all $i$ which implies
$\pi=\|id|$ and $\lambda=\mu$.\qed

If $\lambda$ and $\mu$ are partitions then $\lambda\cle\mu$ is just the
usual inclusion relation among diagrams but in general ``$\cle$'' it is
finer than ``$\subseteq$''.

We proceed by describing the minimal elements lying above $\lambda$.
For a subset $I=\{i_1,\ldots,i_r\}$ of $\{1,\ldots,n\}$ with
$i_1<\ldots<i_r$ we define $c_I(\lambda):=\mu\in\Lambda$ where
$$
\vbox{\halign{$#$\hfill\quad&$#$\hfill\cr
\mu_{i_j}=\lambda_{i_{j+1}}&j=1,\ldots,r-1;\cr
\mu_{i_r}=\lambda_{i_1}+1;&\cr
\mu_i=\lambda_i,&i\not\in I.\cr}}
$$
Clearly, $\lambda\prec c_I(\lambda)$. Conversely,

\Lemma. Let $\lambda,\mu\in\Lambda$ with $\lambda\prec\mu$. Then there
is $I$ such that $c_I(\lambda)\cle\mu$.

\Proof: Let $\pi$ be as in the definition of $\lambda\cle\mu$. If
$\pi=\|id|$ then $\lambda_i<\mu_i$ for some $i$ and we can choose
$I=\{i\}$. Assume from now on $\pi\ne\|id|$.

Let $i_1$ be minimal with $\pi(i_1)\ne i_1$. Then necessarily
$\pi(i_1)>i_1$. Put $i_k:=\pi^{k-1}(i_1)$ and let $r\ge2$ be minimal
with $i_{r+1}>i_r$. Then we have $i_1\le i_r<\ldots<i_2$. Now we take
$I:=\{i_1,i_r,\ldots,i_2\}$ (note that $i_r$ might be equal to $i_1$). For
showing $c_I(\lambda)\cle\mu$ one checks easily that the following
permutation $\pi'$ is defining:
$$
\pi'(i_1)=i_{r+1};\quad\pi'(i_j)=i_j,\quad j=2,\ldots,r;\quad
\pi'(i)=\pi(i),\quad i\not\in I.
$$\qed

\noindent The next lemma shows that one doesn't have to check all
permutations to show $\lambda\cle\mu$.

\Lemma. Let $\lambda,\mu\in\Lambda$ with $\lambda\cle\mu$. Then
$\pi=w_\mu w_\lambda^{-1}$ is a defining permutation.

\Proof: First note that $\pi$ is the permutation with
$k_i(\lambda)=k_{\pi(i)}(\mu)$ for all $i$. Fixing $\pi$ like that
certainly defines a new order relation $\cle'$ which is coarser that
$\cle$. To show that these relations coincide it suffices to show
$\lambda\cle'\mu:=c_I(\lambda)$ for all $I$.

For this we may assume that $I$ is maximal with
$\mu=c_I(\lambda)$. This means that
$$
\eqalign{
&\lambda_i\ne\lambda_{i_1}\phantom{+1}\hbox{ for }i=1,\ldots,i_1-1;\cr
&\lambda_i\ne\lambda_{i_2}\phantom{+1}\hbox{ for }i=i_1+1,\ldots,i_2-1;\cr
&\hbox{etc.}\cr
&\lambda_i\ne\lambda_{i_r}\phantom{+1}\hbox{ for }i=i_{r-1}+1,\ldots,i_r-1;\cr
&\lambda_i\ne\lambda_{i_1}+1\hbox{ for }i=i_r+1,\ldots,n.\cr}
$$
In this case one verifies easily
$$
k_{i_j}(\mu)=k_{i_{j+1}}(\lambda)\quad j=1,\ldots,r-1;\quad
k_{i_r}(\mu)=k_{i_1}(\lambda)
$$
and $k_i(\mu)=k_i(\lambda)$ otherwise. This shows $\lambda\cle'\mu$.
\qed

\Definition: A subset $S\subseteq\Lambda$ is called {\it closed\/} if
$\lambda\in S$ and $\lambda\cle\mu$ implies $\mu\in S$. If that is the
case let $I_S\subseteq\cP$ be the ideal of functions $f$ which vanish
in all $\lamq$ with $\lambda\in\Lambda\setminus S$.

\Theorem Stab. Let $S\subseteq\Lambda$ be a closed subset. Then
$\Xi_i(I_S)\subseteq I_S$ for all $i=1,\ldots,n$.

\Proof: Let $\mu\in\Lambda\setminus S$ and $f\in I_S$. We have to show
that $\Xi_if(\muq)=0$. The definition of $\Xi_i$ shows that
$\Xi_if(\muq)$ is a linear combination of $f(\muq)$ and
$y:=\sigma_i\ldots\sigma_{n-1}\Delta\sigma_1\ldots\sigma_{i-1}f(\muq)$
where each operator $\sigma_j$ is either $s_i$ or $1$. This shows
that $y=f(\lamq)$ with $\mu=c_I(\lambda)$ for some $I$ and
$\lambda\in\Lambda$. Since $S$ is closed we have
$\lambda\in\Lambda\setminus S$ and therefore, $y=0$.\qed

\noindent Now we show the extra vanishing theorem:

\Theorem T1. Let $\lambda,\mu\in\Lambda$ with $\lambda\not\!\!\cle\mu$.
Then $E_\lambda(\muq)=0$.

\Proof: Consider the closed subset
$S=\{\nu\in\Lambda\mid\lambda\cle\nu\}$. We have to show
$E_\lambda\in I_S$. For generic $q$ and $t$ there is a
function $f\in I_S$ with $f(\lamq)\ne0$. Indeed, take for example
$$
f(z):=\prod_{\pi\in W}\left[
\prod_{i<\pi(i)}\phi_{\lambda_i+1}(\lamq_i^{-1}z_{\pi(i)})   
\prod_{i\ge\pi(i)}\phi_{\lambda_i}(q\lamq_i^{-1}z_{\pi(i)})   
\right]
$$
where $\phi_k(z):=(z-1)(z-q^{-1})\ldots(z-q^{-k+1})$. Since $I_S$
is $\Xi_i$\_stable there is $E_{\lambda'}\in I_S$ with
$E_{\lambda'}(\lambda)\ne0$. In particular, $|\lambda'|\le|\lambda|$.
On the other hand, $E_{\lambda'}(\lamq')\ne0$ implies $\lambda'\in S$,
i.e. $\lambda\cle\lambda'$. Therefore, $\lambda'=\lambda$
(\cite{L1}).\qed

\Theorem T2. Let $S\subseteq\Lambda$ be closed. Then
$I_S=\oplus_{\lambda\in S}kE_\lambda$.

\Proof: By \cite{Stab}, we have $I_S=\oplus_{\lambda\in S'}kE_\lambda$
for some subset $S'\subseteq\Lambda$. Let $\lambda\in S'$. Then
$E_\lambda(\lamq)\ne0$ implies $\lambda\in S$, hence $S'\subseteq S$.
Conversely, let $\lambda\in S$ and $\mu\in\Lambda\setminus S$. Then
$E_\lambda(\muq)=0$ by the extra vanishing theorem. Hence
$E_\lambda\in I_S$ and $\lambda\in S'$.\qed

\Corollary. Let $E_\lambda E_\mu=\sum_\nu c_{\lambda\mu}^\nu E_\nu$.
Then $c_{\lambda\mu}^\nu=0$ unless $\lambda,\mu\cle\nu$.

\Proof: Let $S$ be the set of $\nu\in\Lambda$ with $\lambda\cle\nu$.
Then by \cite{T1} and \cite{T2}, the principal ideal $\cP E_\lambda$ is
contained in $\oplus_{\nu\in S}kE_\nu$. This shows $\lambda\cle\nu$
whenever $c_{\lambda\mu}^\nu\ne0$. The relation $\mu\cle\nu$ follows
by symmetry.\qed

\noindent The whole discussion has also a symmetric counterpart. As
already mentioned, on $\Lambda^+$ the order relation $\cle$ is just
inclusion of diagrams. Then everything works for this order relation.
See \cite{KS} for the precise statements.
\def\cPq{{\overline\cP}}

\beginsection Inversion. The inversion formula and integrality results

We have seen that the Macdonald polynomials are obtained from the
$E_\lambda$ or $P_\lambda$ by taking the top homogeneous part. In
this section we show how to invert this process. Let $\Psi:\cP\pf\sim\cP$
be the linear isomorphism mapping $\Eq_\lambda$ to
$E_\lambda$. Another way of describing $\Psi$ is:
let $\cP_d$ be the set of polynomials of degree $d$ which vanish in
$S(n,d-1)$ and $\cPq_d$ the set of homogeneous polynomials of degree
$d$. Then taking the leading term gives an isomorphism
$\cP_d\pfeil\cPq_d$ which maps $E_\lambda$ to $\Eq_\lambda$. Hence,
$\Psi$ is the inverse of this map. This also shows that $\Psi$ is an
$\cH$\_homomorphism. In particular,
$\Psi(\Pq_\lambda)=P_\lambda$. 

Observe that $\Psi$ commutes with $H_i$ and intertwines the action of
$\Phi$ and $\overline\Phi:=z_n\Delta$. Now we study how $\Psi$ behaves
with respect to multiplication by $z_i$. For this we define, for
$i=1,\ldots,n$, the operators $$
Z_i:=t^{n\choose 2}(z_i\Xi_i-1)\Xi_1\ldots\widehat\Xi_i\ldots\Xi_n.
$$
Observe that $S:=t^{n\choose 2}\Xi_1\ldots\Xi_n$ acts like an ``Euler
operator'' on $\cP$, i.e., it acts on $\cP_d$ as scalar multiplication
by $q^d$. We could write $Z_i=(z_i-\Xi_i^{-1})S$ but note
that $\Xi_i^{-1}$ is not a difference operator. Let $M$ be the
operator which acts on $\cPq_d$ by multiplication with $q^{d\choose
2}$.

\Theorem Diag. For all $d\in\NN$ and for $i=1,\ldots,n$ the following
diagram commutes:
$$
\matrix{
\cPq_d          &\Pf{\Psi M}&\cP_{d}          \cr
\downarrow\Rechts{z_i}&        &\downarrow\Rechts{Z_i}\cr
\cPq_{d+1}        &\Pf{\Psi M}&\cP_{d+1}        \cr}
$$
In particular, $Z_i\Psi M=\Psi M z_i$.

\Proof: According to the definition of $\Xi_i$, we have
$$
z_i\Xi_i-1=H_i\ldots H_{n-1}\Phi H_1\ldots H_{i-1}.
$$
Therefore, $Z_i$ maps indeed $\cP_d$ into $\cP_{d+1}$
(\cite{Van}, \cite{Ei}). Moreover, for the
leading term $\overline{Z_i(f)}=q^d z_i\fQ$,
which is equivalent to the commutativity of the diagram.\qed

\noindent Now we can prove the inversion formula:

\Corollary. a) The operators $Z_1,\ldots, Z_n$ commute pairwise.
\Par\noindent
b) For every $f\in\cP$ holds $\Psi M(f)=f(Z_1,\ldots,Z_n)(1)$
(evaluation at the constant function $1$).

\Proof: By \cite{Diag}, we have $\Xi_i=(\Psi M)z_i(\Psi M)^{-1}$. This
implies a). For b) note
$$
\Psi M(f)=\Psi M(f(z_1,\ldots,z_n))(1)=f(Z_1,\ldots,Z_n)(1).
$$\qed

To remove the denominators in the coefficients of $E_\lambda$ we use
the following normalization. Recall, that the {\it diagram} of
$\lambda\in\Lambda$ is the set of points (usually called {\it boxes})
$s=(i,j)\in\ZZ^2$ such that $1\le i\le n$ and $1\le j\le\lambda_i$. For
each box $s$ we define the {\it arm\_length} $a(s)$ and {\it leg-length}
$l(s)$ as
$$
\eqalign{
a(s)&\:=\lambda_i-j\cr
l'(s)&\:=\#\{k=1,\ldots, i-1\mid j\le \lambda_k+1\le\lambda_i\}\cr
l''(s)&\:=\#\{k=i+1,\ldots,n\mid j\le \lambda_k\le\lambda_i\}\cr
l(s)&\:=l'(s)+l''(s)\cr}
$$
If $\lambda\in\Lambda^+$ is a partition then $l'(s)=0$ and
$l''(s)=l(s)$ is just the usual leg\_length. We define
$$ 
\cE_\lambda:=\prod_{s\in\lambda}(1-q^{a(s)+1}t^{l(s)+1})E_\lambda.
$$
$$
\cP_\lambda:=\prod_{s\in\lambda}(1-q^{a(s)}t^{l(s)+1})P_\lambda.
$$
With this normalization, we obtain:
\Proposition. The coefficients of $\cE_\lambda$ and $\cP_\lambda$ are in
$\ZZ[q,q^{-1},t,t^{-1}]$.

\Proof: The leading terms $\overline\cE_\lambda$ and
$\overline\cP_\lambda$  have coefficients in
$\ZZ[q,t]$ by Corollary~5.2 and Theorem~6.1 of
\cite{Kn1}. The result follows from the fact that the operators $Z_i$ are
defined over $\ZZ[q,q^{-1},t,t^{-1}]$.\qed 

\noindent This result can be improved. For $m=1,\ldots,n$ we define the
operators
$$
\eqalign{
A_m&:=H_mH_{m{+}1}\ldots H_{n{-}1}\Phi\cr
\Aq_m&:=\Hq_m\Hq_{m{+}1}\ldots\Hq_{n{-}1}\Phi\cr}
$$
Then we have the following recursion relation:

\Theorem. Let $\lambda\in\Lambda$ with $m:=l(\lambda)>0$. Put
$\lambda^*:=(\lambda_m-1,\lambda_1,\ldots,\lambda_{m{-}1},0,\ldots,0)$.
Then $\cE_\lambda=q^{\lambda_m-1}[\Aq_m-\lamq_mt^m
A_m]\cE_{\lambda^*}$.

\Proof: Apply $\Psi$ to both sides of \cite{Kn1}~Theorem~5.1.\qed

\Corollary. Let $\cE_\lambda=\sum_\mu e_{\lambda\mu}z^\mu$ and
$\cP_\lambda=\sum_\mu p_{\lambda\mu}m_\mu$. Then
$$
e_{\lambda\mu},p_{\lambda\mu}\in t^{-(n-1)(|\lambda|-|\mu|)}\ZZ[q,t].
$$

\Proof: The operators $H_i$ and $\Hq_i$ are defined over $\ZZ[t]$.
Moreover, $\Phi=(z_n-t^{-n+1})\Delta$. Hence the recursion formula
implies $e_{\lambda\mu}\in
t^{-k}\ZZ[q,q^{-1},t]$ with $k=(n-1)(|\lambda|-|\mu|)$. To show that no
negative powers of $q$ appear observe that $z_1\Xi_1-1=H_1\ldots
H_{n-1}\Phi$. Therefore,
$$
q^{\lambda_m-1}\Phi(E_{\lambda^*})=
q^{\lambda_m-1}H_{n-1}^{-1}\ldots H_1^{-1}(z_1(\lamq_1^*)^{-1}-1)
\cE_{\lambda^*}.
$$
The claim follows from $\lamq_1^*=q^{\lambda_m-1}t^{-k}$ for some
$k\in\NN$.

The assertion for $\cP_\lambda$ follows by symmetrization as in the
proof of \cite{Kn1}~Theorem~6.1.\qed

\def\li#1{\mathop{\buildrel#1,r\over\longrightarrow}}
\def\rhoS{{\tilde\rho}}
\def\lamS{{\tilde\lambda}}
\def\muS{{\tilde\mu}}
\def\DeltaS{{\tilde\Delta}}
\def\PhiS{{\tilde\Phi}}
\def\PsiS{{\tilde\Psi}}
\def\XiS{{\tilde\Xi}}
\def\ES{{\tilde E}}
\def\PS{{\tilde P}}
\def\ZS{{\tilde Z}}

\beginsection Classical. The classical limit

Let $\alpha$ be a formal parameter. If one puts $q=t^\alpha$ and lets
$t$ go to $1$ then $(1-t)^{-|\lambda|}\overline\cP_\lambda$ converges to
the Jack polynomial $J^{(\alpha)}_\lambda(z)$. In this final section, I
will discuss the analogue for our non\_homogeneous Macdonald
polynomials. For this it is a bit more convenient to set $t=q^r$ and let
$q$ tend to $1$. Then $\alpha$ and $r$ are related by $\alpha=1/r$.

We introduce the following notation: Let
$p(q,t)\in\QQ(q,t)$, $p_0\in\QQ$ and $k\in\NN$. Then we write $p\li k
p_0$ if $\|lim|\limits_{q\pfeil1}{p(q,q^r)\over(q-1)^k}=p_0$.
For example, $q^at^b-1\li1 a+br$.

As $q\pfeil1$, the points $\lamq$ all collapse to $1$. Therefore, we
introduce the function $\phi_q(x):={x-1\over q-1}$ and the affine
transformation $\phi_q(z):=(\phi_q(z_1),\ldots,\phi_q(z_n))$. Set
$\rhoS:=(0,-r,\ldots-(n-1)r)$. Then we have $\phi_q(\lamq)\li0\lamS$
where $\lamS:=\lambda+w_\lambda\rho$. Conversely, we can write
$\lamq=q^\lamS$.

The action of $\phi_q$ on $\cP$ is given by
$\phi_qf(z)=f(\phi_q^{-1}(z))$. We extend our notation as follows: if
$X(q,t)$ and $X(r)$ are operators on $\cP$ then we write $X(q,t)\li k
X(r)$ if for all $f\in\cP$ $$
\|lim|\limits_{q\pfeil1}(q-1)^{-k}\phi_q X(q,q^r)\phi_q^{-1} f=X(r)f.
$$
For example, if $X(q,t)$ is the multiplication
operator by $g(z;q,t)$ then $\phi_q X(q,q^r)\phi_q^{-1}$ is
multiplication by  $\phi_q g=g\circ\phi_q^{-1}$. In this way, $g\li
k g_0$ has to be understood. Since $\phi_q^{-1}(x)=(q-1)x+1$ we
obtain, for example, $z_i\li01$ while $z_i-1\li1z_i$.

Define the operators $\DeltaS f(z):=f(z_n-1,z_1,\ldots,z_{n-1})$ and
$\PhiS:=(z_n+(n-1)r)\DeltaS$.

\Lemma. We have $\Delta\li0\DeltaS$ and $\Phi\li1\PhiS$.

\Proof: We have
$\phi_q\Delta\phi_q^{-1}f(z)=f(\psi(z_n),z_1,\ldots,z_{n-1})$ where
$\psi(z)=z_n/q-1/q\li0z_n-1$. Moreover,
$z_n(\phi_q^{-1}(z))-t^{-n+1}=(q-1)z_n+(1-t^{-n+1})\li1z_n+(n-1)r$.
\qed

\Theorem. For every $\lambda\in\Lambda$ there is a unique polynomial
$\ES_\lambda$ of degree $|\lambda|$ which vanishes at all $\muS$ with
$|\mu|\le|\lambda|$ and $\mu\ne\lambda$. Moreover,
$E_\lambda\li{|\lambda|}\ES_\lambda$.

\Proof: Repeat the proof of \cite{IntNsym} with respect to vanishing at
$\phi_q(\lamq)$. Then one sees, by induction, that the limit $q\pfeil1$
exists.\qed

\noindent In the limit $q\pfeil1$, the symmetric version has been
already treated in \cite{KS}. As above we obtain

\Theorem. For every $\lambda\in\Lambda^+$ there is a unique polynomial
$\PS_\lambda$ of degree $|\lambda|$ which vanishes at all
$\muS=\mu+\rho$ where $\mu\in\Lambda^+$ with $|\mu|\le|\lambda|$ and
$\mu\ne\lambda$. Moreover,
$P_\lambda\li{|\lambda|}\PS_\lambda$.

\noindent
Next we study the limit of the Hecke operators. Let
$\sigma_i:=s_i-rN_i$.

\Lemma. We have $H_i,\Hq_i\li0\sigma_i$ and $H_i-\Hq_i\li1r$.

\Proof: First observe that $s_i$ commutes with $\phi_q$. Moreover,
$$
(1-t){z_{i+1}\over z_i-z_{i+1}}\phi_q^{-1}=
{1-t\over q-1}{(q-1)z_{i+1}+1\over
z_i-z_{i+1}}\,\li0\,{-r\over z_i-z_{i+1}}.
$$
This implies the claim for $\Hq_i$. For $H_i$, use $H_i=\Hq_i+t-1$.
\qed

\noindent The braid relations for the $H_i$ imply them for the
$\sigma_i$. Moreover, from $H_i\Hq_i=t$ we deduce $\sigma_i^2=1$.
Hence we obtain

\Corollary. The mapping $s_i\mapsto\sigma_i$ extends uniquely to an
action of $W$ on $\cP$.

\Remark: It is not difficult to see that the standard action of $W$ and
the one defined above are conjugated by an element $U\in
\|End|_{\cP^W}\cP\cong M_{n!}(\cP^W)$. It would be interesting to
find such a $U$ explicitly. Note however that $U$ is not uniquely
determined.

\medskip\noindent
Observe that the commutation relations $z_{i+1}H_i=\Hq_i z_i$ can we
rewritten as $(z_{i+1}-1)H_i=\Hq_i(z_i-1)-(H_i-\Hq_i)$ which implies
$$
z_{i+1}\sigma_i=\sigma_i z_i -r;\qquad
z_j\sigma_i=\sigma_i z_j,\quad j\ne i,i+1.
$$
Therefore, the $\sigma_i$ and $z_j$ generate a graded Hecke algebra.

We now consider the limit of the Cherednik operators. Let
$$
\XiS_i:=z_i-\sigma_i\ldots\sigma_{n-1}\PhiS\sigma_1\ldots\sigma_{i-1}.
$$
Then we have

\Theorem. For $i=1,\ldots,n$ holds $\Xi_i-1\li1-\XiS_i$. Moreover,
$\XiS_i\ES_\lambda=\lamS_i\ES_\lambda$ for all $\lambda\in\Lambda$.
In particular, the $\XiS_i$ commute pairwise.

\Proof: Follows from $z_i^{-1}-1\li1-z_i$, $H_i\li0\sigma_i$,
$\Phi\li1\PhiS$, and $\lamq^{-1}-1\li1-\lamS$.\qed

\Corollary. The top homogeneous part of $\PS_\lambda$ and
$\ES_\lambda$ is a Jack polynomial and Opdam's
non\_symmetric analogue {\rm\cite{Op}}, respectively.

\noindent
The extra vanishing theorem goes through verbatim:

\Theorem. Let $\lambda,\mu\in\Lambda$ with $\lambda\not\!\!\cle\mu$.
Then $\ES_\lambda(\muS)=0$.

\noindent
For the inversion formula we introduce
$$
\ZS_i:=z_i-\XiS_i=
\sigma_i\ldots\sigma_{n-1}\PhiS\sigma_1\ldots\sigma_{i-1}.
$$
Furthermore, let $\PsiS:\cP\pfeil\cP$ be the linear automorphism which
maps the leading term of $\ES_\lambda$ to $\ES_\lambda$. Then we
obtain

\Theorem. We have $Z_i\li1\ZS_i$. Moreover, $\ZS_i\Psi=\Psi z_i$ and
the inversion formula holds: $\Psi(f)=f(\ZS_1,\ldots,\ZS_n)(1)$ for all
$f\in\cP$.

\noindent Finally, as for integrality, we choose the following
normalizing factors:
$$
\tilde\cE_\lambda:=
\prod_{s\in\lambda}\big((a(s)+1)+(l(s)+1)r\big)\ES_\lambda.
$$
$$
\tilde\cP_\lambda:=
\prod_{s\in\lambda}\big(a(s)+(l(s)+1)r\big)\PS_\lambda
$$
With this normalization, we obtain:

\Theorem. We have
$\cE_\lambda\li{2|\lambda|}(-1)^{|\lambda|}\tilde\cE_\lambda$ and 
$\cP_\lambda\li{2|\lambda|}(-1)^{|\lambda|}\tilde\cP_\lambda$.
Moreover, the coefficients of $\cE_\lambda$ and $\cP_\lambda$ are in
$\ZZ[r]$.

\beginsection References. References

\baselineskip12pt
\parskip2.5pt plus 1pt
\hyphenation{Hei-del-berg}
\def\L|Abk:#1|Sig:#2|Au:#3|Tit:#4|Zs:#5|Bd:#6|S:#7|J:#8||{%
\edef\TEST{[#2]}
\expandafter\ifx\csname#1\endcsname\TEST\relax\else
\immediate\write16{#1 hat sich geaendert!}\fi
\expandwrite\AUX{\neverexpand\ref{#1}{\TEST}}
\HI{[#2]}
\ifx-#3\relax\else{#3}: \fi
\ifx-#4\relax\else{#4}{\sfcode`.=3000.} \fi
\ifx-#5\relax\else{\it #5\/} \fi
\ifx-#6\relax\else{\bf #6} \fi
\ifx-#8\relax\else({#8})\fi
\ifx-#7\relax\else, {#7}\fi\Par}

\def\B|Abk:#1|Sig:#2|Au:#3|Tit:#4|Reihe:#5|Verlag:#6|Ort:#7|J:#8||{%
\edef\TEST{[#2]}
\expandafter\ifx\csname#1\endcsname\TEST\relax\else
\immediate\write16{#1 hat sich geaendert!}\fi
\expandwrite\AUX{\neverexpand\ref{#1}{\TEST}}
\HI{[#2]}
\ifx-#3\relax\else{#3}: \fi
\ifx-#4\relax\else{#4}{\sfcode`.=3000.} \fi
\ifx-#5\relax\else{(#5)} \fi
\ifx-#7\relax\else{#7:} \fi
\ifx-#6\relax\else{#6}\fi
\ifx-#8\relax\else{ #8}\fi\Par}

\def\Pr|Abk:#1|Sig:#2|Au:#3|Artikel:#4|Titel:#5|Hgr:#6|Reihe:{%
\edef\TEST{[#2]}
\expandafter\ifx\csname#1\endcsname\TEST\relax\else
\immediate\write16{#1 hat sich geaendert!}\fi
\expandwrite\AUX{\neverexpand\ref{#1}{\TEST}}
\HI{[#2]}
\ifx-#3\relax\else{#3}: \fi
\ifx-#4\relax\else{#4}{\sfcode`.=3000.} \fi
\ifx-#5\relax\else{In: \it #5}. \fi
\ifx-#6\relax\else{(#6)} \fi\PrII}
\def\PrII#1|Bd:#2|Verlag:#3|Ort:#4|S:#5|J:#6||{%
\ifx-#1\relax\else{#1} \fi
\ifx-#2\relax\else{\bf #2}, \fi
\ifx-#4\relax\else{#4:} \fi
\ifx-#3\relax\else{#3} \fi
\ifx-#6\relax\else{#6}\fi
\ifx-#5\relax\else{, #5}\fi\Par}
\setHI{[WUN]$\,$}
\sfcode`.=1000

\L|Abk:Ch|Sig:Ch|Au:Cherednik, I.|Tit:Non\_symmetric
Macdonald's polynomials|Zs:Preprint|Bd:-|S:25 pages|J:1995||

\L|Abk:HU|Sig:HU|Au:Howe, R., Umeda, T.|Tit:The Capelli
identity, the double commutant theorem, and multiplicity\_free
actions|Zs:Math. Ann.|Bd:290|S:569--619|J:1991||

\L|Abk:Kn1|Sig:Kn|Au:Knop, F.|Tit:Integrality of two
variable Kostka functions|Zs:Preprint|Bd:-|S:12 pages|J:1996||

\L|Abk:KS|Sig:KS|Au:Knop, F.; Sahi, S.|Tit:Difference
operators and symmetric functions defined by their zeros%
|Zs:Preprint|Bd:-|S:14 pages|J:1996||

\B|Abk:M1|Sig:M1|Au:Macdonald, I.|Tit:Symmetric functions and Hall
polynomials (2nd ed.)|Reihe:-|Verlag:Clarendon Press|Ort:Oxford|J:1995||

\L|Abk:M2|Sig:M2|Au:Macdonald, I.|Tit:Affine Hecke algebras and
orthogonal polynomials|Zs:S\'eminaire Bourbaki|Bd:-|S:n$^o$ 797|J:1995||

\L|Abk:Ok|Sig:Ok|Au:Okounkov, A.|Tit:Quantum immanants and
higher Capelli identities|Zs:Preprint|Bd:-|S:23 pages|J:1995||

\L|Abk:Ol|Sig:Ol|Au:Olshanski, G.|Tit:Quasi\_symmetric functions and
factorial Schur functions|Zs:Preprint|Bd:-|S:21 pages|J:1995||

\L|Abk:Op|Sig:Op|Au:Opdam, E.|Tit:Harmonic analysis for certain
representations of graded Hecke algebras%
|Zs:Acta Math.|Bd:175|S:75--121|J:1995||

\Pr|Abk:Sa|Sig:Sa|Au:Sahi, S.|Artikel:The spectrum of certain invariant
differential operators associated to Hermitian symmetric spaces%
|Titel:Lie theory and geometry|Hgr:J.-L. Brylinski et al. eds.%
|Reihe:Progress Math.|Bd:123|Verlag:Birkh\"auser|Ort:Boston%
|S:569--576|J:1994||

\L|Abk:Wa|Sig:WUN|Au:Wakayama, M.; Umeda, T.; Noumi, M.|Tit:A quantum
analogue of the Capelli identity and an elementary differential calculus
on $GL_q(n)$|Zs:Duke Math. J.|Bd:76|S:567--594|J:1994||

\bye